\documentstyle[12pt,amscd,amssymb]{amsart}

\newcommand{\seq}[3]{{#1}_{#2}, \ldots, {#1}_{#3}}

\newcommand{\Dws}{D^{(w,\Sigma)}}
\newcommand{\M}{{\cal N}_g} 
\newcommand{\Y}{\Sigma \times {{\Bbb S}}^1}
\newcommand{\X}{X^o}
\newcommand{\Spz}{\text{Sp}\, (2g,{\Bbb Z})}
\newcommand{\Cabg}{\CC[\a,\b,\g]}

\newcommand{\surj}{\twoheadrightarrow}
\newcommand{\inc}{\hookrightarrow}
\newcommand{\ar}{\rightarrow}
\newcommand{\bd}{\partial}
\newcommand{\x}{\times}
\newcommand{\ox}{\otimes}
\newcommand{\iso}{\cong}
\newcommand{\isom}{\stackrel{\simeq}{\ar}}

\newcommand{\CP}{{\Bbb C \Bbb P}}

\newcommand{\End}{\text{End}}
\newcommand{\Diff}{\text{Diff}}

\newcommand{\Sym}{\text{Sym}}

\newcommand{\PD}{\text{P.D.}}

\newcommand{\cN}{{\cal N}}

\newcommand{\cS}{{\cal S}}

\newcommand{\cU}{{\cal U}}

\renewcommand{\AA}{{\Bbb A}}
\newcommand{\CC}{{\Bbb C}}
\newcommand{\DD}{{\Bbb D}}

\newcommand{\TT}{{\Bbb T}}
\newcommand{\SS}{{\Bbb S}}
\newcommand{\ZZ}{{\Bbb Z}}

\renewcommand{\a}{\alpha}
\renewcommand{\b}{\beta}

\newcommand{\g}{\gamma}

\renewcommand{\l}{\lambda}

\newcommand{\s}{\sigma}

\newcommand{\p}{\psi}
\newcommand{\q}{\phi}

\newcommand{\z}{\zeta}

\renewcommand{\S}{\Sigma}
\newcommand{\D}{\Delta}
\renewcommand{\L}{\Lambda}
\renewcommand{\P}{\Phi}

\newcommand{\frg}{{\frak g}}

\newcommand{\frm}{{\frak m}}

\theoremstyle{plain}
\begingroup
\theorembodyfont{\sl}
\newtheorem{thm}{Theorem}
\newtheorem{cor}[thm]{Corollary}
\newtheorem{lem}[thm]{Lemma}
\newtheorem{prop}[thm]{Proposition}
\newtheorem{conj}[thm]{Conjecture}
\endgroup

\theoremstyle{definition}

\theoremstyle{remark}
\newtheorem{rem}[thm]{Remark}
\newtheorem{ex}[thm]{Example}
\newtheorem{notation}[thm]{Notation}

\title{Ring structure of the Floer cohomology of $\Y$}

\author{Vicente Mu\~noz}
\address{Departamento de \'Algegra, Geometr\'{\i}a y Topolog\'{\i}a \\ Facultad de Ciencias \\ Universidad de M\'alaga \\ 29071 M\'alaga \\ Spain}
\email{vmunoz@@agt.cie.uma.es}
\thanks{\hbox{$^*$}Supported by a grant from Ministerio de Educaci\'on y Cultura 
}
\date{September, 1997}

\begin{document}

\maketitle

\begin{abstract}
  We give a presentation for the Floer cohomology ring $HF^*(\Y)$, where $\S$ is a
  Riemann surface of genus $g \geq 1$, which coincides with the conjectural 
  presentation for the quantum cohomology ring of the moduli space of flat
  $SO(3)$-connections of odd degree over $\S$. We study the spectrum of the action of
  $H_*(\S)$ on $HF^*(\Y)$ and prove a physical assumption made in~\cite{Vafa}.
\end{abstract}

\section{Introduction}
\label{sec:intro}

Let $\S=\S_g$ be a Riemann surface of genus $g \geq 1$ and let $\M$ denote the moduli
space
of flat $SO(3)$-connections with nontrivial second Stiefel-Whitney class $w_2$. 
This is a smooth symplectic manifold
of dimension $6g-6$. Alternatively, we can consider $\S$ as a smooth complex curve of
genus $g$ and $\M$ as the moduli space of odd degree rank two 
stable vector bundles on $\S$ with fixed determinant, which is a smooth
complex variety of complex dimension $3g-3$. The symplectic deformation class of $\M$ 
only depends on the genus $g$ and not on the particular complex structure on $\S$.

We consider the following rings associated to the Riemann surface
(we will always use $\CC$-coefficients):

\begin{itemize}
\item $QH^*(\M)$ is the quantum cohomology of $\M$ (see~\cite{RT}). This
is well-defined since $\M$ is a
positive symplectic manifold. As vector spaces, $QH^*(\M)=H^*(\M)$,
but the multiplicative structure is different. The minimal 
Chern number of $\M$ is $2$, so $QH^*(\M)$ is $\ZZ/4\ZZ$-graded (the grading 
comes from reducing mod $4$ the $\ZZ$-grading of $H^*(\M)$). The ring structure of 
$QH^*(\M)$, called quantum multiplication, 
is a deformation of the usual cup product for $H^*(\M)$. It is associative
and graded commutative. We remark that we do not introduce Novikov rings to define
$QH^*(\M)$ as in~\cite[section 8]{RT}
(otherwise said, as $H^2(\M) =\ZZ$, we should introduce an extra variable $q$ of degree $4$, then we equate $q=1$).

\item $HF_*^{\mathrm{symp}}(\M)$ is the symplectic Floer homology of $\M$ (with the
symplectomorphism $\q=\text{id})$. The symplectic 
manifold $\M$ is connected, simply connected and $\pi_2(\M)=\ZZ$
 (see~\cite[introduction]{DS}), so the groups
$HF_*^{\mathrm{symp}}(\M)$ are well-defined~\cite{Floer1}. They are $\ZZ/4\ZZ$-graded.
$HF_*^{\mathrm{symp}}(\M)$ is endowed with the pair of pants 
product~\cite{Piunikhin}, which is an associative and graded commutative ring
structure. The symplectic Floer cohomology of $\M$, $HF^*_{\mathrm{symp}}(\M)$,
is defined as the dual of the symplectic Floer homology. There is a Poincar\'e 
duality~\cite[remark 2.4]{Piunikhin} and a pairing $<,>$.

\item $HF_*(\Y)$ is the instanton Floer homology of the three manifold $Y=\Y$ for the $SO(3)$-bundle with second Stiefel-Whitney class
$w_2=\PD[\SS^1] \in H^2(Y;\ZZ/2\ZZ)$. This is defined in~\cite{Floer2} and is
$\ZZ/4\ZZ$-graded. We introduce a multiplication on $HF_*(Y)$
using a suitable four-dimensional cobordism~\cite[section 5]{Salamon}. 
Let $X$ be the four manifold given as
a pair of pants times $\S$, which is
a cobordism between $Y \sqcup \, Y$ and $-Y$. This gives a map 
$HF_*(Y) \ox HF_*(Y) \ar HF_*(Y)$, which
is an associative and graded commutative ring structure on $HF_*(Y)$.
Again, the instanton Floer cohomology of $Y$, $HF^*(Y)$, is the dual of $HF_*(Y)$.
There is a pairing $<,>:HF^*(Y) \ox HF^*(-Y) \ar \CC$. As $Y=\Y$ admits a orientation
reversing self-diffeomorphism, we can identify $HF^*(-Y) \iso HF^*(Y)$, and hence
we have a pairing on $HF^*(Y)$ (see~\cite{Don1}~\cite{Don2}).
We will denote $HF_g^*=HF^*(\Y)$.
\end{itemize}

\begin{thm}
\label{thm:mainiso}
There are natural isomorphisms of vector spaces
\begin{equation} 
   QH^*(\M) \iso HF^*_{\mathrm{symp}}(\M) \iso HF^*(\Y).
\label{eqn:mainiso}
\end{equation}
Moreover the first isomorphism respects the ring structures.
\end{thm}

\begin{pf}
The second isomorphism is due to Dostoglou and Salamon~\cite[theorem 10.1]{DS}.
It is the particular case 
where one considers $\q=\hbox{id}:\S \ar \S$, in which the mapping torus of $\q$ is
$\Y$ and the $SO(3)$-bundle has $w_2=\PD [\SS^1]$. 
The second isomorphism is a standard result
obtained by Floer~\cite{Floer1}.
In~\cite[theorem 5.1]{Piunikhin} it is proved that the first
isomorphism intertwines the products. 
\end{pf}

\begin{conj}
\label{conj:iso}
  The second isomorphism in~\eqref{eqn:mainiso} is a ring isomorphism.
\end{conj}

In~\cite{Salamon} D. Salamon announced a proof of conjecture~\ref{conj:iso}
but finally he could not complete it. Many implications of this conjectural result
to four-dimensional topology were given by Donaldson~\cite{Don2}. The author followed
this program in~\cite{Thesis} for small genus. Later he obtained very
nice results on the behaviour of Donaldson invariants under the operation
of connected sum along a Riemann surface~\cite{genus2}~\cite{genusg}, 
exploiting only the isomorphism~\eqref{eqn:mainiso} as vector spaces.

Using physical methods, Vafa et al.\ ~\cite{Vafa} find a set of 
generators and relations for $QH^*(\M)$. There are two main assumptions in their
argument. The first one is conjecture~\ref{conj:iso}. The second
one is that the spectrum of the action of $H_*(\S)$
on $HF^*(\Y)$ can be read off from the Donaldson invariants of $\S \x \TT^2$.

Later Siebert and Tian~\cite{Siebert} claimed to have found a mathematical proof
for the presentation of $QH^*(\M)$ 
given in~\cite{Vafa} but they could not yet finish their program. 
In this paper we prove that the set of generators and relations given in~\cite{Vafa}
is a presentation for $HF^*(\Y)$ (theorem~\ref{thm:main}), 
following a method inspired in that of Siebert and Tian~\cite{Siebert}. 
This together with completion of the work~\cite{Siebert} will produce a
proof of conjecture~\ref{conj:iso} (although in a rather indirect way).
We also prove the physical assumption on the spectrum of $HF^*(\Y)$ in~\cite{Vafa}
(proposition~\ref{prop:15}).

We leave the implications of theorem~\ref{thm:main} to Donaldson invariants of four-manifolds (mostly in the case $b^+=1$) for future work.

{\em Acknowledgements.\/} I am very grateful to Bernd 
Siebert and Gang Tian for providing me with a copy of~\cite{Siebert} which was very enlightening. In particular, the proof of theorem~\ref{thm:8}
is due entirely to them.
Thanks to the organization of the CIME Course on 
Quantum Cohomology held in Cetraro (Italy, 1997) for inviting me.
Also discussions with Paul Seidel were useful.
Finally I acknowledge the hospitatility of the Mathematics 
Department in Universidad de M\'alaga.

\section{Ring structure of $H^*(\M)$}
\label{sec:2}

Let us recall the known description of the homology of $\M$~\cite{King}~\cite{ST2}.
Let $\cU \ar \S \x \M$ be the universal bundle and consider the K\"unneth
decomposition as in~\cite{King}
$$
   c_2(\End_0 \, \cU)=2a [\S] + 4 \p -b
$$
with $\p=\sum c_i \g_i^{\#}$, where  $\{\seq{\g}{1}{2g}\}$ is a symplectic basis of 
$H_1(\S;\ZZ)$ with $\g_i \g_{i+g}=[\S]$ for $1 \leq i \leq g$,
and $\{\g_i^{\#}\}$ is the dual basis of $H^1(\S)$. 
In terms of the map $\mu: H_*(\S) \to H^{4-*}(\M)$,
given by $\mu(a)= -{1 \over 4} \,
p_1(\frg_{\,\cU}) / 4$ (here $\frg_{\,\cU} \to \S \x \M$ is the associated
universal $SO(3)$-bundle, and $p_1(\frg_{\,\cU}) \in H^4(\S\x\M)$ its first
Pontrjagin class), we have
$$
   \left\{ \begin{array}{l} a= 2\, \mu(\S) \in H^2
   \\ c_i= \mu (\g_i) \in H^3, \qquad 0\leq i \leq 2g
   \\ b= - 4 \, \mu(x) \in H^4       
    \end{array} \right.
$$
where $x \in H_0(\S)$ is the class of the point, and $H^i=H^i(\M)$.
These elements generate $H^*(\M)$ as a 
ring~\cite{King}~\cite{Thaddeus}. 
So there is a basis $\{f_s\}_{s \in \cS}$ for $H^*(\M)$ with elements of the form
$$
   f_s=a^nb^mc_{i_1}\cdots c_{i_r},
$$
for a finite set $\cS$ of multi-indices of the
form $s=(n,m; i_1,\ldots,i_r)$, $n,m \geq 0$, $r \geq 0$, $1 \leq
i_1 < \cdots < i_r \leq 2g$. 
The mapping class group $\Diff(\S)$ acts on $H^*(\M)$, with the action factoring
through the action of $\Spz$ on $\{c_i\}$. 
The invariant part, $H_I^*(\M)$, is generated by
$a$, $b$ and $c=-2 \sum_{i=0}^g c_ic_{i+g}$. Then
\begin{equation}
   H_I^*(\M)= \CC [a, b, c]/I_g,
\label{eqn:k}
\end{equation}
where $I_g$ is the ideal of relations satisfied by $a$, $b$ and $c$. 
Here $\deg (a)=2$, $\deg(b)=4$, $\deg(c)=6$. 
Actually, a basis for $H_I^*(\M)$ is given by the monomials $\a^a\b^b\g^c$, $a+b+c
<g$~\cite{ST2}. For
$0 \leq k \leq g$, the primitive component of $\L^k H^3$ is 
$$
   \L_0^k H^3 = \ker (c^{g-k+1} : \L^k H^3 \ar \L^{2g-k+2} H^3).
$$
Then the $\Spz$-decomposition of $H^*(\M)$ is~\cite{King} 
$$
   H^*(\M)= \bigoplus_{k=0}^g \L_0^k H^3 \ox \CC [a, b, c]/I_{g-k}.
$$
\begin{prop}[\cite{King}\cite{ST2}]
\label{prop:3}
  For $g=1$, let $q_1^1=a$, $q_1^2=b$, $q_1^3=c$. Define recursively, for $g \geq 1$,
$$
   \left\{ \begin{array}{l} q_{g+1}^1= a q_g^1 +g^2 q_g^2
   \\ q_{g+1}^2 = b q_g^1 + {2g \over g+1}  q_g^3
   \\ q_{g+1}^3 = c q_g^1
    \end{array} \right.
$$
  Then $I_g=(q_g^1, q_g^2, q_g^3)$, for all $g \geq 1$.
\end{prop}

\begin{pf}
  Define $\z_0=1$ and
  $\z_{n+1}=a \z_n +n^2 b \z_{n-1} + 2n(n-1) c \z_{n-2}$, for all $n \geq 0$. 
  Theorem 3.1 in~\cite{King} establishes that
  $I_g=(\z_g, \z_{g+1}, \z_{g+2})$. The proposition follows from the obvious equalities
$$
   \left\{ \begin{array}{l}  q_g^1=\z_g
   \\  q_g^2= {1 \over g^2}(\z_{g+1} - a \z_g)
   \\  q_g^3= {1 \over 2g(g+1)} \left( \z_{g+2}-a \z_{g+1} -(g+1)^2 b\z_g \right))
    \end{array} \right.
$$  
\end{pf}

\section{Ring structure of $HF^*_g$}
\label{sec:3}

With the aid of the basis $\{f_s \}_{s \in \cS}$ for $H^*(\M)$ we are going 
to construct a basis for $HF^*_g=HF^*(\Y)$ to understand its ring structure.
We need to use the gluing properties of the Floer homology of a three manifold.
Put $Y=\Y$ and let $w_2= \PD[\SS^1] \in H^2(Y;\ZZ/2\ZZ)$.
Let us state the result that we shall use.

\begin{prop}[\cite{Don1}~\cite{Don2}~\cite{Thesis}]
\label{prop:4}
  For any smooth oriented four-manifold $\X$ 
  with boundary $\bd \X=Y$, any $w \in H^2(X;\ZZ)$ with 
  $w|_Y = \PD[\SS^1]$ in $H^2(Y;\ZZ/2\ZZ)$, and any $z \in \AA(\X)$, we have defined
  a relative invariant $\q^{w}(\X, z ) \in HF_*(Y)$. These relative invariants enjoy
  the following gluing property, suppose $X=\X_1 \cup_Y \X_2$ is a 
  closed four-manifold split into two open four manifolds $\X_i$ with $\bd \X_1= Y$,  
  $\bd \X_2 =-Y$, and $w \in H^2(X;\ZZ)$ satisfying $w|_Y =\PD[\SS^1]$ in 
  $H^2(Y;\ZZ/2\ZZ)$. Put $w_i =w|_{X_i}$. Then for
  $z_i \in \AA(\X_i)$, $i=1,2$, we have
  \begin{equation}
     D_X^{(w,\S)}(z_1z_2)=
    <\q^{w_1}(\X_1, z_1 ),\q^{w_2}(\X_2, z_2 )>, \label{eqn:sym}
  \end{equation}
  where  $D^{(w,\S)}_X=D^w_X +D^{w+\S}_X$ ($D_X^w$ is
  the Donaldson invariant of $X$ for $w$, see also~\cite{Thesis}~\cite{genus2}).
  When $b^+=1$, the invariants are calculated for a long neck, i.e. we
  refer to the invariants defined by $\S$.
\end{prop}

Consider the manifold $A=\S \x D^2$,
with boundary $Y=\Y$, and let $\D= \text{pt} \x D^2 \subset A$ be the 
horizontal slice with $\bd \D=\SS^1$. Put $w=\PD [\D] 
\in H^2(A;\ZZ)$. Clearly $\AA(A)= \AA(\S)= \Sym^*(H_0(\S)
\oplus H_2(\S)) \otimes \bigwedge^* H_1(\S)$. For every 
$s \in \cS$, $f_s=a^n b^mc_{i_1}\cdots c_{i_r}$, define
\begin{eqnarray*}
   z_s&=& \S^n x^m \g_{i_1} \cdots \g_{i_r} \in \AA(\S), \\
   e_s&=& \q^w (A, z_s ) \in HF^*(Y)=HF^*_g
\end{eqnarray*}
(here we identify Floer homology and Floer cohomology through Poincar\'e duality).
Then $\{e_s \}_{s \in \cS}$ is a basis for $HF^*_g$. This is a consequence 
of~\cite[lemma 21]{genus2}. The product $HF_g^* \ox HF_g^* \ar HF_g^*$ 
is given by $\q^w(A, z_S) \q^w(A,z_{S'})=\q^w(A,z_Sz_{S'})$. 
Then $\q^w(A,1)$ defines the neutral element of the product.
As a consequence, the following elements are generators of $HF_g^*$,
\begin{equation}
   \left\{ \begin{array}{l} \a= 2 \, \q^w(A,\S) \in HF^2_g
   \\ \p_i= \q^w(A,\g_i) \in HF^3_g, \qquad 0\leq i \leq 2g
   \\ \b= - 4 \, \q^w(A,x) \in HF^4_g
    \end{array} \right.
\label{eqn:j}
\end{equation}
Note that there is an obvious epimorphism of rings $\AA(\S) \surj HF_g^*$.

\begin{thm}
\label{thm:deform}
  Denote by  $*$ the product induced in $H^*(\M)$ by the product in $HF_g^*$ 
  under the isomorphism $H^*(\M) \isom HF_g^*$ given by $f_s\mapsto e_s$, 
  $s \in \cS$. Then $*$ is a deformation of the cup-product
  graded modulo $4$, i.e. for $f_1 \in H^i (\M)$, $f_2 \in H^j(\M)$, it is
  $f_1 * f_2 = \sum_{r \geq 0} \P_r(f_1,f_2)$, where $\P_r \in H^{i+j-4r}(\M)$ and 
  $\P_0= f_1 \cup f_2$.
\end{thm}

\begin{pf}
  First, for $s, s' \in \cS$,
  $$
   <e_s, e_{s'}> = \Dws_{\S \x \CP^1}(z_sz_{s'}) =0,
  $$
  unless $\deg(f_s) + \deg(f_{s'})= 6g-6 +4r$, $r \geq 0$, as these are the
  only 
  possible
  dimensions for the moduli spaces of anti-self-dual connections on $\S \x
  \CP^1$.
  Moreover, when $\deg(f_s) + \deg(f_{s'})= 6g-6$, the moduli space is
  $\M$, so $<e_s, e_{s'}>=- <f_sf_{s'}, [\M]>= -<f_s, f_{s'}>$ (the
  minus sign is due to the different convention orientation for 
  Donaldson invariants).

  Now let $f_s$, $f_{s'}$ be basic elements of degrees $i$ and $j$ respectively.
  Put $f_sf_{s'}= \sum c_tf_t$ and
  $f_s * f_{s'}= \sum d_tf_t$. This means that $e_se_{s'}= \sum d_te_t$.
  Write $e_se_{s'}= \sum_m g_m$, where $g_m =\sum_{\deg (f_t)=m} d_t e_t$
  are the homogeneous parts. Put $\hat g_m=\sum_{\deg (f_t)=m} d_t f_t$. Let $M$ be
  the maximum $m$ such that $g_m \neq 0$. 
  Then there is $f_{s''}$ of degree $6g-6-M$ such that $<\hat g_M, f_{s''}> \neq 0$.
  Since
  $$  
  0 \neq - <\hat g_M, f_{s''}>= <g_M, e_{s''}>= <e_se_{s'}, e_{s''}>= \Dws_{\S \x
   \CP^1}(z_sz_{s'}z_{s''}),
  $$
  it is $\deg(f_s) + \deg(f_{s'})+ \deg(f_{s''}) \geq 6g-6$, i.e. $M \leq 
  i+j$. Now for $m= i+j$, any $f_{s''}$
  of degree $6g-6-m$, it is $<\hat g_m, f_{s''}>= -\Dws_{\S \x
   \CP^1}(z_sz_{s'}z_{s''})=<f_sf_{s'}f_{s''}, [\M]>=<f_sf_{s'}, f_{s''}>$. So 
  $\hat g_{i+j}=f_sf_{s'}$. 

  Finally, $<e_se_{s'}, e_{s''}>= \Dws_{\S \x
  \CP^1}(z_sz_{s'}z_{s''})=0$, whenever 
  $\deg(f_s) + \deg(f_{s'})+ \deg(f_{s''}) \not\equiv 6g-6 \pmod 4$, so $\hat g_m =0$ 
  unless $m \equiv i+j \pmod 4$.
\end{pf}

\begin{rem}
\label{rem:deform}
  We do not claim that the isomorphism in theorem~\ref{thm:deform} is the one
  in~\eqref{eqn:mainiso}. Actually this is not the case (see example~\ref{ex:17}).
\end{rem}

There is again an action of $\Diff(\S)$ on $HF_g^*$ factoring through an action
of $\Spz$ on $\{\p_i\}$. The invariant part $(HF_g^*)_I$ is 
generated by $\a$, $\b$ and $\g= -2  \sum_{i=0}^g \q^w(A, \g_i \g_{i+g})$.
The epimorphism $\Cabg \surj (HF_g^*)_I$, $z \mapsto \q^w(A,z)$, allows us 
to write
\begin{equation}
    HF^*(\Y)_I=\Cabg /J_g,
\label{eqn:l}
\end{equation}
where $J_g$ is the ideal of relations of $\a$, $\b$ and $\g$.
Now $\deg (\a)=2$, $\deg(\b)=4$, $\deg(\g)=6$, but $J_g$ is not a homogeneous ideal.

\begin{lem}
\label{lem:decom}
Suppose $\g J_g \subset J_{g+1}$, for all $g \geq 1$. Then we have
the $\Spz$-decomposition
$$
   HF^*(\Y)= \bigoplus_{k=0}^g \L_0^k H^3 \ox \CC [\a, \b, \g]/J_{g-k}.
$$
\end{lem}

\begin{pf}
  The isomorphisms in theorem~\ref{thm:mainiso}
  respect the $\Spz$-action and hence induce isomorphisms on the invariant parts.
  Then $\dim (HF^*_g)_I = \dim H^*_I(\M)$, for all $g \geq 1$. Now the lemma 
  is a consequence of the argument in the proof of~\cite[proposition 2.2]{King}
  and the discussion preceding it.
\end{pf}

\section{A presentation for $(HF^*_g)_I$}
\label{sec:4}

Theorem~\ref{thm:deform} and the arguments in~\cite[section 2]{ST}
imply that we can deform the relations of $H^*_I(\M)$ to get a presentation
for $(HF_g^*)_I$. More explicitly, 

\begin{lem}
\label{lem:deform}
  It is $(HF_g^*)_I= \Cabg/ (R^1_g , R^2_g, R^3_g)$, where $R^i_g= q^i_g +$ lower order
  terms of degrees $\deg q_g^i-4r$, $r>0$, 
  as polynomials in $\Cabg$ ($q^i_g$ are defined in proposition~\ref{prop:3}).
\end{lem}

\begin{pf}
  Suppose first that $g \geq 2$.
  Granted theorem~\ref{thm:deform}, \cite[theorem 2.2]{ST} implies that $J_g= (R^1_g ,
  R^2_g, R^3_g)$, where $R_g^i$ is
  $q_g^i$ expressed in terms of $\a$, $\b$ and $\g$ and the multiplication of $HF_g^*$.
  Now we note that under the isomorphism $H^*(\M) \isom HF_g^*$ of
  theorem~\ref{thm:deform}, $a \mapsto \a$, $b \mapsto \b$, $c \mapsto \g$
  (it always can be arranged so that these elements are in the basis, as $g \geq 2$, 
  see~\cite[proposition 4.2]{ST2}). 
  So $R_g^i$ is equal to $q_g^i$ plus lower order terms. 
  The case $g=1$ is computed directly in lemma~\ref{lem:9}.
\end{pf}

\begin{lem}
\label{lem:7}
  $J_{g+1} \subset J_g$, for all $g \geq 1$.
\end{lem}

\begin{pf}
  Let $\S_g$ be a Riemann surface of genus $g$ and consider 
  $$
    \S_{g+1} \subset A_g= \S_g \x D^2 \subset S= \S_g \x \CP^1,
  $$ 
  where $\S_{g+1}$ is given by $\S_g$ with a trivial handle added internally. 
  Then the map $H_*(\S_{g+1}) \ar H_*(\S_g)$ induces $\AA(\S_{g+1}) \ar   
  \AA(\S_g)$ which sends $(\a, \b, \g) \mapsto (\a, \b, \g)$.
  Put $A_{g+1}=\S_{g+1} \x D^2 \subset A_g$. This gives
  a map $(HF^*_{g+1})_I \ar (HF^*_g)_I$,
  $\q^w(A_{g+1}, z) \mapsto \q^w(A_g, z)$. Put $S= S^o \cup_{\S_{g+1} \x \SS^1}
  A_{g+1}$. Let
  $z \in J_{g+1}$. Then $\q^w(\S_{g+1} \x D^2, z)=0$. 
  So for any $z_s \in \AA(\S_g)$, $s \in \cS$, 
$$
  \Dws_S(z \,z_s)= <\q^w(\S_{g+1} \x D^2, z), \q^w(S^o, z_s)> =0.
$$
  This is equivalent to $z \in J_g$.
\end{pf}

\begin{thm}
\label{thm:8}
  There are numbers $c_{g+1}, d_{g+1} \in \CC$ such that, for all $g \geq 1$, 
  $$
   \left\{ \begin{array}{l} R_{g+1}^1 = \a R_g^1 + g^2 R_g^2
   \\ R_{g+1}^2 = (\b+c_{g+1}) R_g^1 + {2g \over g+1}  R_g^3
   \\ R_{g+1}^3 = \g  R_g^1 + d_{g+1} R_g^2
    \end{array} \right.
  $$
\end{thm}

\begin{pf}
  We follow almost literally the argument of 
  Siebert and Tian~\cite[proposition 3.2]{Siebert}.
  As $R_{g+1}^1 \in J_{g+1} \subset J_g$ is a relation on degree $2g+2$, it is a linear
  combination of $\a R_g^1$ and $R_g^2$. Looking at the leading terms 
  (proposition~\ref{prop:3}), we
  have $R_{g+1}^1 = \a R_g^1 + g^2 R_g^2$.
  Analogously, $R_{g+1}^2$ is a combination of $\a^2 R_g^1$, $\b R_g^1$,
  $\a R_g^2$ and $R_g^1$. Only the term $R_g^1$ has degree less than $2g+4$,
  so $R_{g+1}^2 = \b R_g^1 + {2g \over g+1}  R_g^3 + c_{g+1} R_g^1$, for an 
  unknown coefficient $c_{g+1}$.
  In the same fashion, $R_{g+1}^3$ is $\g  q_g^1$  plus a linear
  combination of  $R_g^2$ and $\a R_g^1$. Adding a suitable multiple of
  $R_{g+1}^1$ (which is always allowed without loss of generality), 
  we have $R_{g+1}^3 = \g  q_g^1 + d_{g+1} R_g^2$.
\end{pf}

\begin{lem}
\label{lem:9}
  The starting relations (for $g=1$) are $R_1^1=\a$, $R_1^2=\b-8$ 
  and $R_1^3=\g$.
\end{lem}

\begin{pf}
  $HF^*_1$ is of dimension $1$, i.e. $HF^*_1=\CC$ (see~\cite{Don2}~\cite{Thesis}). 
  Let $S$ be the $K3$ surface and fix
  an elliptic fibration for $S$, whose fibre is be $\S=\TT^2$. The 
  Donaldson invariants are, for $w \in H^2(S; \ZZ)$ with $w \cdot \S \equiv 1
  \pmod 2$ (see~\cite{Thesis}),
  $$
   \Dws_S (e^{tD})=-e^{-Q(tD)/2}.
  $$
  Then $\Dws_S (1)=-1$ and $\Dws_S (\S^d)=0$, for $d >0$. 
  Also from~\cite[remark 4]{genusg}, $\Dws_S (x)=2$.
  Pet $S^o$ be the complement of an open tubular neighbourhood of $\S$ in $S$.
  Then $\q^w(S^o, 1)$ generates $HF_1^*$ and $\q^w(S^o, \S)=0$, 
  $\q^w(S^o, x) =-2 \, \q^w(S^o,1)$ and $\q^w(S^o, \g_1\g_2)=0$,
  i.e. $\a=0$, $\b-8=0$ and $\g=0$ in $HF_1^*$ (recall~\eqref{eqn:j}). 
\end{pf}

\begin{prop}
\label{prop:10}
  For $g \geq 2$, there exists a non-zero vector $v \in HF^*_g$ such that
  \begin{eqnarray*}
    \a v &=& \left\{ \begin{array}{ll} 4(g-1) \, v & \text{$g$ even} \\
                      4(g-1) \sqrt{-1} \, v \qquad & \text{$g$ odd} \end{array}
    \right. \\
    \b v &=& (-1)^{g-1} 8 \, v \\
    \g v & =& 0
  \end{eqnarray*}
\end{prop}

\begin{pf}
  We shall construct such a vector as the relative invariants of an open 
  four-manifold $\X$ with boundary $\bd \X=Y=\Y$, where
  the closed four-manifold $X=\X \cup_Y A$ is of simple type with $b^+>1$ and $b_1=0$.
  For concreteness, let $X$ be the manifold $C_g$ from~\cite[definition 25]{genusg}. 
  We recall its construction. Let $S_g$ denote the elliptic surface of geometric genus
  $p_g=g-1$ and with no multiple fibres. It contains a section  $\s$
  which is a rational curve of self-intersection
  $-g$. Let $F$ be the elliptic fibre. Then $\s +gF$ can be represented by an embedded
  Riemann surface $\tilde{\S}$ of genus $g$ and self-intersection $g$. Blow-up 
  $S_g$ at $g$ points in $\tilde{\S}$
  to get $B_g$ with an embedded Riemann surface $\S_g$ of genus $g$ and
  self-intersection zero. Then put $X=C_g= B_g \#_{\S_g} B_g$ (the double of $B_g$ 
  along $\S_g$).
  By~\cite[proposition 27]{genusg},
  $X$ is of simple type and $\DD^w_X(e^{\a})= D_X^w((1+{x \over 2})e^{\a})=
  -2^{3g-5} e^{Q(\a)/2} e^{K \cdot \a} +
   (-1)^g 2^{3g-5} e^{Q(\a)/2} e^{-K \cdot \a}$, 
  where $K \in H^2(X;\ZZ)$ satisfies $K \cdot \S_g=2g-2$ 
  ($w\in H^2(C_g; \ZZ)$ is a particular element, which we do not need to specify here).
  Let us suppose from now on that $g$ is even, the other case being similar.
  By~\cite[proposition 3]{genusg},
$$
  \Dws_{X}(e^{\a})= -2^{3g-5} e^{Q(\a)/2} e^{K \cdot \a} +
   (-1)^g 2^{3g-5} e^{Q(\a)/2} e^{-K \cdot \a}.
$$
  We set $v = \q^w(\X,\S + 2g-2) \in HF^*(\Y)=HF_g^*$.
  Let us prove that this is the required element. For any $z_s=\S^n x^m 
  \g_{i_1} \cdots \g_{i_r}$, it is~\cite[remark 4]{genusg},
$$
  <v, e_s>=  \Dws_{X}((\S+2g-2) z_s)= 
   \left\{ \begin{array}{ll} 0, & r >0 \\ -2^{3g-4}(2g-2)^{n+1} 2^m,\qquad & r=0
   \end{array} \right.
$$
  Then $<\a v , e_s> = <\q^w(\X,2 \,\S(\S + 2g-2)), \q^w(A, z_s)>= 
  \Dws_X((\S+2g-2) 2 \,\S\, z_s)= (4g-4) <v, e_s>$, for all $s \in 
  \cS$. Then $\a v = (4g-4) v$. Analogously,
  $\g v =0$ and $\b v= -8 v$.
\end{pf}

\begin{notation}
\label{not:ini}
  We set $R^1_0=1$, $R^2_0=0$ and $R^3_0=0$.
\end{notation}

\begin{thm}
\label{thm:11}
  For all $g \geq 1$, $c_g=(-1)^g 8$ and $d_g=0$.
\end{thm}

\begin{pf}
  The result is true for $g=1$ by lemma~\ref{lem:9} and notation~\ref{not:ini}.
  Suppose it is true for $1 \leq r \leq g$, and let us prove it for $g+1$. By 
  proposition~\ref{prop:10}, there exists $v \in HF_{g+1}^*$ with
  $\b v=(-1)^g 8 \, v$, $\g v = 0$ and $\a v = 4g \, v$ if $g$ is odd and $\a 
  v= 4g \sqrt{-1} \, v$ if $g$ is even.

  In first place, $\g v=0$ implies $R_r^3 v=0$, for $1 \leq r \leq g$.
  In second place, $\b v=(-1)^g 8 \, v$ implies
\begin{eqnarray*}
  R_g^2 v &=& (\b +(-1)^g 8) R_{g-1}^1 v =(-1)^g 16 R_{g-1}^1 v, \\
  R_{g-1}^2 v &=& (\b +(-1)^{g-1} 8) R_{g-1}^1 v =0, \\
  R_{g-2}^2 v &=& (-1)^g 16 R_{g-3}^1 v, \\ R_{g-4}^2 v &=& 0, \\ \vdots
\end{eqnarray*}
  In third place,
  $R_g^1 v = \a R_{g-1}^1 v+ (g-1)^2 R_{g-1}^2 v = \a R_{g-1}^1 v$,
  $R_{g-2}^1 v =\a R_{g-3}^1 v$, $\ldots$
  Also 
  $$
   R_{g-1}^1 v = \a R_{g-2}^1 v+ (g-2)^2 R_{g-2}^2 v = (\a^2 +(g-2)^2
   (-1)^g 16 )R_{g-3}^1 v.
  $$
  So finally,
  $$ 
  R_{g-1}^1 v = \left\{
  \begin{array}{ll}
  (\a^2 + (-1)^g 16 (g-2)^2) \cdots (\a^2 +(-1)^g 16 \cdot 1^2 ) v & 
  \text{$g$ odd} \\
  (\a^2 + (-1)^g 16 (g-2)^2) \cdots (\a^2 +(-1)^g 16 \cdot 2^2 ) \a v \quad& 
  \text{$g$ even}
  \end{array} \right.
  $$
  As a conclusion $R_{g-1}^1 v = \l v$, with $\l \neq 0$, and
  $$ 
  \left\{
  \begin{array}{l} R_g^1 v = \a R_{g-1}^1 v  \\ R_g^2 v = (-1)^g 16
  R_{g-1}^1 v \\ R_g^3 v =0 
  \end{array} \right.
  $$
  As $v \in HF_{g+1}^*$, we have $R_{g+1}^1 v =0$, $R_{g+1}^2 v =0$ and 
  $R_{g+1}^3 v =0$. Evaluate the equations from theorem~\ref{thm:8} on $v$ to get
  $c_{g+1}=(-1)^{g+1} 8$ and $d_{g+1}=0$.
\end{pf}

\begin{cor}
\label{cor:12}
  We have $\g J_g \subset J_{g+1} \subset J_g$, for all $g \geq 1$.
\end{cor}

\begin{pf}
  The second inclusion is lemma~\ref{lem:7}. For the first inclusion, note that
  $\g R^1_g =R^3_{g+1} \in J_{g+1}$ by the third equation in theorem~\ref{thm:8}. 
  Then multiplying the first two equations in theorem~\ref{thm:8} we get that
  $\g R^2_g, \g R^3_{g} \in J_{g+1}$.
\end{pf}

  Using this corollary in lemma~\ref{lem:decom}, we have finally proved that
\begin{thm}
\label{thm:main}
  The Floer cohomology of $\Y$, for $\S=\S_g$ a Riemann surface of genus $g$, has
  a presentation
  $$
   HF^*(\Y)= \bigoplus_{k=0}^g \L_0^k H^3 \ox \Cabg /J_{g-k}.
  $$
  where $J_r=(R^1_r, R^2_r,R^3_r)$ and $R^i_r$ are defined recursively by setting 
  $R^1_0=1$, $R^2_0=0$, $R^3_0=0$ and putting for all $r \geq 0$
  $$
   \left\{ \begin{array}{l} R_{r+1}^1 = \a R_r^1 + r^2 R_r^2
   \\ R_{r+1}^2 = (\b+(-1)^{r+1}8) R_r^1 + {2r \over r+1}  R_r^3
   \\ R_{r+1}^3 = \g  R_r^1
    \end{array} \right.
  $$
\end{thm}

\begin{rem}
\label{rem:14}
  The presentation obtained for $HF_g^*$ is the conjectural presentation for
  $QH^*(\M)$ (see~\cite{Siebert}).
\end{rem}

\begin{cor}
\label{cor:13}
  $\ker (\g:(HF_g^*)_I \ar (HF_g^*)_I ) = J_{g-1}/J_g \subset \Cabg /J_g 
  =(HF_g^*)_I$.
\end{cor}

\begin{pf}
  By the corollary~\ref{cor:12}, $\g$ factors as
  $$ 
    \Cabg/J_g \surj \Cabg /J_{g-1} \stackrel{\g}{\inc} \Cabg/J_g.
  $$
  The second map is a monomorphism since $\a^a\b^b\g^c$, $a+b+c<g-1$, form a basis
  for $\Cabg /J_{g-1}$, and their image under $\g$ are linearly independent in
  $\Cabg /J_g$. 
  The corollary follows.
\end{pf}

  For any $F \in \Cabg$ define the expectation value by $<F>_g=<F_g,
  1>_{HF_g^*}$, where $1 \in HF_g^*$ is the unit element. Therefore
  $<F_1, F_2>_{HF_g^*}=<F_1 F_2>_g$.

\begin{cor}
\label{cor:13bis}
  For any $F \in \Cabg$, $<\g F>_g= -2g <F>_{g-1}$.
\end{cor}

\begin{pf}
  By corollary~\ref{cor:12}, the formula above holds for any $F \in J_{g-1}$,
  as both sides are zero. So it is enough to check it for a set of elements
  generating $HF_{g-1}^*$, i.e. for $F_{abc}=
  \a^a\b^b\g^c$, $a+b+c<g-1$. If $(a,b,c)
  \neq (0,0,0)$, it is $<F_{abc}>_{g-1}=0$ and $<\g F_{abc}>_g=0$ by 
  degree reasons. Now $<\g^{g-1}>_{g-1}=-<c^{g-1}, [{\cN}_{g-1}]>$ hence
  the corollary follows from $<c^g, [\M]>=-2g <c^{g-1}, [{\cN}_{g-1}]>$
  (see~\cite{Thaddeus}).
\end{pf}

\section{Local ring decomposition of $HF^*(\Y)$}
\label{sec:local}

In~\cite{Vafa} it is asserted that the only eigenvalues of the action of 
$\mu(\S)$, $\mu(x)$ and $\mu(\g_i)$ on $HF^*(\Y)$ are the ones given by looking at the Donaldson invariants of the manifold $X=\S \x \TT^2$ i.e.
if we denote by $W \subset HF_g^*$ the image $\q^w(\X,\AA(\S))$, where
$X=\X \cup_Y A$, then $\a$, $\b$ and $\g$ act on $W$ and their eigenvalues
are all the eigenvalues of their action on $HF_g^*$. 
The following result is
a proof of this physical assertion.

\begin{prop}
\label{prop:15}
  The eigenvalues of $(\a, \b,\g)$ in $(HF_g^*)_I$ are
  $(0,8,0)$, $( \pm 4, -8 ,0)$, $(\pm 8 \sqrt{-1}, 8,0)$, $\ldots$, 
  $(\pm 4(g-1) \sqrt{-1}^g, (-1)^{g-1} 8 , 0)$.
\end{prop}

\begin{pf}
  Put $V=(HF_g^*)_I$.
  As $\g J_{g-1} \subset J_g$, one has $\g^g \in J_g$, i.e. $\g^g=0$ in
  $V$, so the only eigenvalue of $\g$ is zero. To compute the eigenvalues of
  $\a$, $\b$ we can restrict to $V/ \g V$
  (if $p$ is a polynomial with $p(\a)=0$ in $V/\g V$, 
  then $p(\a)$ is a multiple of $\g$ in
  $V$ and $p(\a)^g=0$ in $V$).
  Now the ideal of relations of $V$ can be written as $J_g=(\z_g,\z_{g+1}, \z_{g+2})$,
  where $\z_0=1$ and
  $\z_{r+1}=\a \z_r +r^2 (\b+(-1)^r 8) \z_{r-1} + 2r(r-1) \g \z_{r-2}$, 
  for all $r \geq 0$ 
  (see proposition~\ref{prop:3}). So 
  $$  
    V/\g V= \CC[\a, \b]/(\bar{\z}_g , \bar{\z}_{g+1}),
  $$
  where $\bar{\z}_0=1$, $\bar{\z}_{r+1}=\a \bar{\z}_r + r^2 (\b +(-1)^r8)  
  \bar{\z}_{r-1}$, for
  $r \geq 0$. From $r^2 (\b +(-1)^r8) \bar{\z}_{r-1}= \bar{\z}_{r+1} - \a \bar{\z}_r$
  we infer that $(\b +(-1)^r8) \bar{\z}_{r-1} \in (\bar{\z}_r , \bar{\z}_{r+1})$. 
  Continuing in this way, 
  $$
   (\b +(-1)^g8)(\b +(-1)^{g-1}8) \cdots (\b  - 8) \in (\bar{\z}_g , \bar{\z}_{g+1}),
  $$
  which implies that the only eigenvalues of $\b$ in $V/\g V$, and hence in $V$,
  are $\pm 8$.
  Let us study the eigenvalues of $\a$ for $\b=8$, $\g=0$. Again we
  only need to study
  $V/(\g, \b-8)V= \CC[\a]/(\hat{\z}_g , \hat{\z}_{g+1})$, where now
  $\hat{\z}_0=1$, $\hat{\z}_{r+1}=\a \hat{\z}_r + r^2 (8+(-1)^r8) \hat{\z}_{r-1}$.
  Then
  $$
   \left\{ \begin{array}{ll} 
   \hat{\z}_r=(\a^2 +(r-2)^216) \cdots (\a^2 + 2^2 16)\a^2 \qquad & \text{$r$ even} \\
   \hat{\z}_r=(\a^2 +(r-1)^216) \cdots (\a^2 + 2^2 16)\a & \text{$r$ odd} 
   \end{array} \right.
  $$
  from where the eigenvalues of $\a$ will be $0, \pm 8\sqrt{-1}, \pm 16\sqrt{-1}, 
  \ldots, \pm 8 \left[{g-1 \over 2}\right]\sqrt{-1}$.
  We leave the other case to the reader.
\end{pf}

\begin{rem}
\label{rem:16}
  As mentioned in~\cite{Vafa}, by the very definition of $\g=-2 \sum \q^w(A, \g_i
  \g_{i+g})$, it is $\g^{g+1}=0$ in $HF_g^*$, so the only eigenvalue of $\g$ is zero.
\end{rem}

Proposition~\ref{prop:15} says that $(HF_g^*)_I$ can be decomposed as a sum of local 
artinians rings
\begin{equation}
  (HF_g^*)_I =\bigoplus_{i=-(g-1)}^{g-1} R_{g,i},
\label{eqn:x}\end{equation}
where $R_{g,i}$ is a local artinian ring with maximal ideal
$\frm=(\a-4i,\b - (-1)^i 8,\g)$ if $i$ is odd, 
$\frm=(\a-4i \sqrt{-1},\b - (-1)^i 8,\g)$ if $i$ is even.
Also $HF_g^*$ is decomposed as 
\begin{equation}
  HF_g^* =\bigoplus_{k=0}^g \bigoplus_{i=-(g-k-1)}^{g-k-1} \L_0^kH^3 \ox R_{g-k,i} =
  \bigoplus_{i=-(g-1)}^{g-1} \bigoplus_{k=0}^{g-|i|-1} \L_0^kH^3 \ox R_{g-k,i}.
\label{eqn:y}\end{equation}

We recall from lemma~\ref{lem:9} that $HF_1^*= \Cabg/(\a, \b-8, \g)$. Let us see the
next cases.

\begin{ex}
\label{ex:17}
  For $g=2$, 
  $J_2=(\a^2+\b-8 , \a(\b +8) + \g , \a\g)$. In $(HF_2^*)_I$,
  $\g= -\a(\b +8)$ and $\g\a=0$ yield $\a^2 (\b+8)=0$. Now $\a^2=-(\b-8)$ so
  $(\b-8)(\b+8)=0$ and
  $\a^2(\a^2-16)=0$. Also $\g J_1 \subset J_2$ implies 
  $\g\a=\g^2=\g (\b-8)=0$. Finally $(\g+16\a)(\a^2-16)=-16\g +16\a(\a^2-16)= 
  -16(\g+\a(\b+8))=0$. All together proves
$$ 
  (HF_2^*)_I= {\Cabg \over (\a-4, \b+8, \g)} \oplus {\Cabg \over (\a^2, \b- 8,
    \g+16\a)} \oplus {\Cabg \over (\a+4, \b+8, \g)}.
$$
  We want to remark that $HF_2^* \isom QH^*(\cN_2)$
  (see~\cite[example 5.3]{Piunikhin} for a presentation of the
  latter ring). The isomorphism sends $\a \mapsto h_2$, $\b \mapsto -4(h_4-1)$,
  $\g \mapsto 4(h_6-h_2)$, where $h_2$, $h_4$, $h_6$ are the generators
  of $QH^2$, $QH^4$, $QH^6$ respectively. This 
  was conjectured in~\cite[conjecture 1.22]{Thesis}.
\end{ex}

\begin{ex}
\label{ex:18}
  For $g=3$, 
  $J_3=\big( \a (\a^2+\b-8) + 4(\a\b+8\a +\g) , (\b-8)(\a^2 + \b -8) 
  +{4 \over 3}\a \g
  , \g (\a^2 + \b-8)\big)$. Put $V=(HF_3^*)_I$. Then 
  $$ 
   V/\g V= \CC[\a,\b]/(\a (\a^2+\b-8) + 4(\a\b+8\a ) , (\b-8)(\a^2 + \b -8)).
  $$
  In $V/\g V$, the first relation yields $-5 \a (\b-8)= \a^3 + 64\a $
  and the second
  $\a (\b- 8)(\a^2 + \b -8) = 0$. This implies $\a (\a^2-16 )(\a^2+64)=0$.
  Also $(\b-8) \a (\a^2-16)=0$. Using $(\b-8)\a^2 =-(\b-8)^2$, we get
  $\a (\b-8)(\b+8)=0$.

  Therefore, in $V$, 
  $\a (\a^2-16 )(\a^2+64)$  and $\a (\b-8)(\b+8)$
  are multiples of $\g$. As $\g J_2
  \subset J_3$, we have $\g \a^2 (\a^2 -16)=0$ and $\g\a^2(\b-8)=0$
  by example~\ref{ex:17}. So
  $\a^3 (\a^2-16 )^2(\a^2+64)=0$, $\a^3 (\a^2 -16) (\b-8)(\b+8)=0$ and
  $\a^3 (\b-8 )^2(\b+8)=0$. It can be checked now that
$$ 
  (HF_3^*)_I= { \Cabg \over (\a- 8 \sqrt{-1}, \b-8, \g)} \oplus
   {\Cabg \over ((\a-4)^2,
   \b+ 8, \g+8(\a-4))} \oplus 
$$
$$ \oplus {\Cabg\over (\a^3, \a(\b-8), (\b-8)^2-{64 \over 3}\a^2, \g+16\a)}
    \oplus {\Cabg \over ((\a+4)^2,
   \b+ 8, \g+8(\a+4))} \oplus  
$$
$$
  \oplus {\Cabg \over (\a+ 8 \sqrt{-1}, \b-8, \g)}.
$$
\end{ex}

\section{Conjecture}
\label{sec:conj}

We state the following conjecture, which first occurred to Paul Seidel and
the author in mid'96.

\begin{conj}
\label{conj:ultima}
The decomposition in equation~\eqref{eqn:y} is
$$
 HF_g^* \iso H^*( s^0 \S) \oplus H^* (s^1\S) \oplus \cdots \oplus
 H^* (s^{g-2}\S) \oplus 
$$
$$
 \oplus H^* (s^{g-1}\S) \oplus H^* (s^{g-2}\S) \oplus 
 \cdots \oplus H^*( s^0 \S),
$$
where $s^i\S$ is the $i$-th symmetric product of $\S$.
Here $H^*(s^i\S)$ is isomorphic to the eigenspace of eigenvalues
$(\pm 4(g-1-i) \sqrt{-1}^{g-i},(-1)^{g-1-i}8,0)$. The isomorphism respects only
$\ZZ/2\ZZ$-grading and is $\Diff(\S)$-equivariant.
\end{conj}

Simple computations establish that the dimensions of both vector spaces appearing
in conjecture~\ref{conj:ultima} are the
same, i.e. $g2^g$. The Euler characteristic are both vanishing. Moreover the 
dimensions of the invariant parts coincide ${g+2 \choose 3}$. 
Examples~\ref{ex:17} and~\ref{ex:18} agree with the conjecture.

A deeper reason for the above conjecture is the fact that $HF_g^*$ is the
space for a gluing theory of Donaldson invariants associated to 
the three manifold $Y=\Y$. The gluing theory of Seiberg-Witten invariants
should be based on the Seiberg-Witten-Floer homology groups of $\Y$, which
are indexed by a line bundle $L$ (the determinant line bundle of the 
$\text{spin}^c$-structure on $Y$). The only possibilities are $c_1(L)=\pm 
(g-1-i)[\SS^1]$, $0 \leq i \leq g-1$ (see~\cite[section 6]{Thesis}). It is
believed that the Seiberg-Witten-Floer groups for $L$ are isomorphic to 
$H^*(s^i\S)$.

\end{document}